\documentclass[10pt]{article}
\usepackage{amsmath}
\usepackage{amssymb}
\usepackage{enumerate}
\usepackage{theorem}
\usepackage{amsfonts}
\usepackage[a4paper,margin=2cm,vmargin={1cm,2cm},includeheadfoot]{geometry}
\usepackage{ifthen}
\usepackage[scaled=0.92]{helvet}
\usepackage[latin1]{inputenc}
\usepackage{srcltx}
\usepackage{tikz}
\usepackage[labelfont=bf,margin=1cm,justification=justified,labelsep=period]{caption}
\definecolor{gris}{gray}{0.8}

\newcommand{\buil}[3]{\mathrel{\mathop{\kern0pt#1}\limits_{#2}^{#3}}}

\theoremheaderfont{\bfseries}
\newtheorem{theorem}{Theorem}
\newtheorem{definition}[theorem]{Definition}
\newtheorem{proposition}[theorem]{Proposition}
\newtheorem{lemma}[theorem]{Lemma}

\newtheorem{corollary}[theorem]{Corollary}
\newtheorem{example}[theorem]{Example} 
{\theorembodyfont{\rmfamily}
\newtheorem{remark}[theorem]{Remark}
}

\numberwithin{equation}{section} \numberwithin{theorem}{section}
\begin{document}

\title{Stability for a class of semilinear fractional stochastic  integral equations}
\date{}

\author{Allan Fiel\thanks{Corresponding author. Partially supported by the CONACyT fellowship 259100. }\\
Depto. de Control Autom\'{a}tico\\
Cinvestav-IPN\\
Apartado Postal 14-740\\ 07000 M\'{e}xico D.F., Mexico\\ afiel@ctrl.cinvestav.mx \and
 Jorge A. Le\'on\thanks{Partially supported by the CONACyT grant 220303.}
\\ Depto. de Control Autom\'{a}tico\\
Cinvestav-IPN\\
Apartado Postal 14-740\\ 07000 M\'{e}xico D.F., Mexico\\ jleon@ctrl.cinvestav.mx
 \and David M\'arquez-Carreras\thanks{Partially supported by the MTM2012-31192 ``Din\'amicas
Aleatorias'' \textit{del Ministerio de Econom\'ia y competitividad}.}\\ Departament
de Probabilitat, L\`{o}gica i Estad\'{\i}stica\\ Universitat de
Barcelona, Gran Via 585\\ 08007-Barcelona, Catalunya \\ davidmarquez@ub.edu} 
\maketitle

\begin{abstract} 

 {\small In this paper we study some stability criteria for some semilinear  integral equations with a function 
 as initial condition and with additive noise, which is a Young integral that could be a functional of fractional Brownian
 motion. Namely, we consider stability in the mean, asymptotic stability, stability, global stability and Mittag-Leffler
 stability. To do so, we use comparison results for fractional equations and an equation (in terms of Mittag-Leffler 
 functions) whose family of solutions includes those of the underlying equation.}
\end{abstract}

\section{Introduction}

Currently fractional systems are of great interest because of the applications they have in several areas of science and technology, such as
engineering, physics, chemistry, mechanics, etc. (see, e.g. \cite{Da}, \cite{H}, \cite{KST}, \cite{Pod1} and the references therein). 
Particularly we can mention system identification \cite{Da}, robotics \cite{MA}, control \cite{Da, Pod1}, electromagnetic theory \cite{Ho},
chaotic dynamics and synchronization \cite{GG, HLQK, YL, YLS}, applications on viscoelasticity \cite{Ba}, analysis of electrode processes
\cite{INK}, Lorenz systems \cite{GG}, systems with retards \cite{DLC}, quantic evolution of complex systems \cite{KBD}, numerical methods
for fractional partial differential equations \cite{DEOK, MO, MO2}, among other. A nice survey of basic properties of deterministic fractional
differential equations is in Lakshmikantham and Vatsala \cite{LV}. Also, many researchers have established stability criteria of mild solutions 
of stochastic fractional differential equations using different techniques.

For deterministic systems, the stability of fractional linear equations has been analyzed by Matignon \cite{M} and Radwan et al. \cite{RSES}.
Besides, several authors have studied non-linear cases using  Lyapunov method (see, e.g. Li et al. \cite{LCP} and its references). In particular,
non-linear  fractional systems with a function as initial condition using also the Lyapunov technique have been considered in the Ph.D.Thesis of
Mart\'{\i}nez-Mart\'{i}nez \cite{MM-L-FA}. Moreover, in the work of Junsheng et al. \cite{JJM} the form of the solution for a linear fractional
equation with a constant initial condition in terms of Mittag-Leffler function is given by means of the Adomian decomposition method. Wen et al.
\cite{WWL} have established stability results for fractional non-linear equations via the Gronwall inequality. Lemma \ref{solueq} below can be
seen as an extension of the results in \cite{YLS} and the Gronwall inequality stated in \cite{WWL}. In \cite{WWL}, the stability is
used to obtain synchronization of fractional systems. 

On the other hand, a process used frequently in literature is fractional Brownian motion $B^H=\{ B_t^H, t\ge 0\}$ due to the wide range
of properties it has, such as long range memory (when the Hurst  parameter $H$ is greater than one half) and intermitency
(when $H<1/2$). Unfortunately, in general, it is not a semimartingale (the exception is $H=1/2$). Thus, we cannot use classical It\^o calculus  
in order to integrate processes with respect to $B^H$ when $H\neq 1/2,$ but we may use another approaches such as Young integration (see
Gubinelli \cite{G}, Young \cite{Y}, Z\"ahle \cite{Z}, Dudley and Norvai\u{s}a \cite{DNo}, Lyons \cite{Ly}). The reader can also see  Nualart
\cite{N1}, and  Russo and Vallois \cite{RV} for other types of integrals. As a consequence, an important application is the analysis of 
stochastic integral equations driven by fractional Brownian motion that has been considered by several authors these days for different 
interpretations of stochastic integrals (see, e.g. Lyons \cite{Ly},  Quer-Sardanyons and Tindel \cite{QT}, Le\'on and Tindel \cite{LT}, 
Nualart \cite{N1}, Friz and Hairer \cite{FH}, Lin \cite{Lin} and Nualart and R\u{a}\c{s}canu \cite{NR}). 

Stability of stochastic systems driven by Brownian motion has been also studied. Some authors use fundamental solution of this equations in 
order to investigate the stability of random systems. An example of this is the paper of Applebay and Freeman \cite{AF}, who give the solution 
in terms of the principal matrix of  integrodifferential equations with an It\^o integral noise and find the equivalence between almost sure
exponential convergence and the $p$-th mean exponential convergence to zero for these systems. Bao \cite{B} uses Gronwall inequality to state
the mean square stability for Volterra-It\^o equations with a function as initial condition and bounded kernels.
Several researchers have studied stability of stochastic
systems via Lyapunov function techniques. An example of this is the paper of  Li et al. \cite{LLW}, who proves stability in probability for  
It\^o-Volterra integral equation, also Zhang and Li \cite{ZL} have stated a stochastic type stability criteria for stochastic integrodifferential
equations with infinite retard, and, Zhang and Zhang \cite{ZZ} have dealed with conditional stability of Skorohod Volterra type equations with
anticipative kernel. Nguyen \cite{DN} present the solution via the fundamental solution for linear stochastic differential equations with 
time-varying delays to obtain the exponential stability of these systems. The noise is an additive one and has form $\int_0^\cdot 
\sigma(s)dW^H_s.$  Here  
$$W_t^H=\int_0^t (t-s)^{H-1/2}dW_s,\quad H\in (1/2,1),$$
$W$ is a Brownian motion and $\sigma$ is a deterministic function such that
$$\int_0^\infty \sigma^2(s)e^{2\lambda s}ds<\infty ,$$
for some $\lambda>0.$ Also, Zeng et al. \cite{ZYC} utilize the Lyapunov function techniques to prove  stability in probability and moment
exponential stability for stochastic differential equation driven by fractional Brownian motion with parameter $H>1/2$ .
Yan and Zhang \cite{YZ} proved sufficient conditions for the asymptotical stability in $p$-th moment for the closed form of the solution to a
fractional impulsive partial neutral stochastic integro-differential equation with state dependent retard in Hilbert space. In the linear case, Fiel 
et al. \cite{FLM} have used the Adomian decomposition method to find the mild solution of an stochastic fractional integral equation with a
function as initial condition driven by a H\"older continuous process in terms of Young or Skorohod integrals. This closed form is given in terms 
of Mittag-Leffler functions. The stability in the large and stability in the mean sense of these random systems is also analized. As an application,
the stability of  equations driven by a functional of fractional Brownian motion is derived. 

In this paper we extend the results given in \cite{FLM} and \cite{WWL}, that is, we study the stability of the solution to the  equation
\begin{equation}\label{maineq}
X(t)=\xi_t+\frac{1}{\Gamma(\beta)} \int_0^t (t-s)^{\beta-1} [A X(s)+h(X(s))] ds + 
Z_t,\quad t\ge 0.
\end{equation}
The initial condition  $\xi=\{\xi_t, t\ge 0\}$ is a function, $h$ is a $O(x)$ as $x\to 0$ (i.e. there are $C>0$ and $\delta>0$
such that $|h(x)|\le C|x|$ for $|x|<\delta$), $\beta\in(0,1)$, $A<0$ and
 $Z$ is a Young integral of the form
$$Z_t=\frac{1}{\Gamma(\alpha)}\int_0^t (t-s)^{\alpha-1}d\theta_s.$$
Here, 
$\theta=\{\theta_s, s\ge 0\}$ is a $\gamma$-H\"older continuous function
that  may represent the 
pahts of a functional of fractional Brownian motion, where
 $\gamma\in (0,1)$, $\alpha\in(1,2)$ and $\alpha+\gamma>2$. Unlike other papers where the involved kernels  
are bounded functions we consider the case that kernels are not bounded, and use comparison results as a main
tool. It is worth mentioning that it is considered the stability of the solution to (\ref{maineq}) in \cite{WWL} with
$Z\equiv 0$ and $\xi$ a constant.

This work is organized as follows. In Section \ref{preli} we introduce a fractional integral equation, whose family of solutions
include those of (\ref{maineq}). Also, in Section \ref{preli}, we state a comparison result for fractional systems that 
becomes the main tool for our results. In Section \ref{clasemi1}, we study some stability criteria for equation (\ref{maineq})
in the case that $Z\equiv 0.$ These results can be seen as extensions of the results given in \cite{FLM} and \cite{WWL}.  
Finally, the stability of equation (\ref{maineq}) in the case that $\theta$ is either a H\"older continuous process, or a
functional of fractional Brownian motion is considered in Section \ref{adino1}.


\section{Preliminaries}\label{preli}

In this section we introduce the framework and the definitions that we use to prove our results. Although some 
results are well-known, we give them here for the convenience of the reader. Part of the main tool that we need is
the stability of some fractional linear systems as it was presented by Fiel et al. \cite{FLM} and a comparison result
(see Lemma \ref{pach} below). 

\subsection{The Mittag-Leffler function}
The Mittag-Leffler function is an important tool of fractional calculus due to its properties and applications. 
As we can see in Lemma \ref{solueq} below, solutions of semilinear fractional integral equations depend on it. 
In order to see a more detailed exposition on this function, the reader is refered to the book of Podlubny \cite{Pod1}. 
For $z\in\mathbb{R},$ this function is defined as 
\begin{displaymath}
E_{a,b}(z)=\sum_{k=0}^{\infty}\frac{z^k}{\Gamma(ka +b)},\quad a, b>0,
\end{displaymath}
where $\Gamma$ is the Gamma function. A function $f$ defined on $(0,\infty)$ is said to be completely monotonic
if it possesses derivatives $f^{(n)}$ for all $n\in\mathbb{N}\cup\{ 0\},$ and if
\[ (-1)^nf^{(n)}(t)\ge 0 \]
for all $t>0.$ In particular, we have that each completely monotonic function on $(0,\infty)$ is positive, decreasing 
and convex, with concave first derivative (see, e.g. \cite{MS} and \cite{S}).

It is well-known that for $a,b\ge 0,$ $z\mapsto E_{a,b}(-z)$ is completely monotonic if and only if $a\in (0,1]$ 
and $b\geq a$ (see Schneider \cite{S}). Moreover, for $a<2,$ there is a positive constant $C_{a,b}$ such that 

\begin{equation}\label{estiml}
\Big{|}E_{a,b}(z) \Big{|}\leq \frac{C_{a,b}}{1+|z|},\quad z\le 0,
\end{equation}
(see \cite{Pod1}, Theorem 1.6). For $a,b>0$ and $\lambda\in\mathbb{R},$ this function satisfies the following
(see (1.83) in \cite{Pod1}):
\begin{equation}\label{derimit}
\frac{d}{dz}\left(z^{b-1}E_{a,b}(\lambda z^a)\right)=z^{b-2}E_{a,b-1}(\lambda z^a).
\end{equation}
Also, for $b>0$ (see equality (1.99) in \cite{Pod1}), we have
\begin{equation}\label{intmitlef}
\int_0^z s^{b-1}E_{a,b}(\lambda s^a)ds=z^bE_{a,b+1}(\lambda z^a).
\end{equation}

\subsection{The Young integral}

Here we introduce the Young integral, which is an integral with respect to H\"older continuous functions.
This was initially defined for functions with $p$-variation in Young \cite{Y}.

For $T> 0$ and $\gamma\in (0,1),$ let $C^\gamma_1([0,T];\mathbb{R})$ be the set of 
$\gamma$-H\"older continuous functions 
$g:[0,T]\to\mathbb{R}$ of one variable such that, the seminorm 
$$||g||_{\gamma,[0,T]}:=\sup_{r,t\in [0,T], r\neq t}\frac{|g_t-g_r|}{|t-r|^\gamma},$$
is finite. Also by $||g||_{\infty, [0,T]},$ we denote the supremum norm of $g.$ 

Using the following result, we can understand easily the basic properties of Young integral for H\"older 
continuous functions. The proof of this theorem can be found in \cite{G} (see also \cite{LT}). Sometimes 
we write ${\mathcal J }_{st}(fdg)$ instead of
$\int_s^t f_udg_u.$
\begin{theorem}\label{gubiyoung} Let $f\in C_1^\kappa ([0,T];
\mathbb{R})$ and $g\in C_1^\gamma ([0,T];\mathbb{R}),$ 
with $\kappa+\gamma >1.$ Then:

\begin{itemize}
\item[1.]  ${\mathcal J}_{st}(fdg)$ 
coincides with the usual Riemann integral if $f$ and $g$ are smooth functions.

\item[2.] We have, for $s\le t\le T,$
$${\mathcal J}_{st}(fdg)=\lim_{|\Pi_{st}|\rightarrow 0}\sum_{i=0
}^{n-1}f_{t_i} (g_{t_{i+1}}-g_{t_i}),$$
where the limit is over any partition $\Pi_{st}= \{ t_0=s,\ldots,t_n=t \}$ of $[s,t],$ whose
mesh tends to zero.  
In particular, ${\mathcal J}_{st}(fdg)$ coincides with the Young integral as
defined in \cite{Y}.

\item[3.] The integral ${\mathcal J}(fdg)$ satisfies:
$$|{\mathcal J}_{st}(fdg)|\leq ||f||_{\infty}||g||_\gamma|t-s|^\gamma
+c_{\gamma,\kappa}||f ||_{\kappa}||g||_\gamma |t-s|^{\gamma+\kappa},$$
where $c_{\gamma,\kappa}=(2^{\gamma+\kappa}-2)^{-1}.$

\end{itemize}

\end{theorem}
We observe that this integral has been extended by Z\"ahle \cite{Z}, Gubinelli \cite{G}, 
Lyons \cite{Ly}, among others. For a more detailed exposition on the Young integral the 
reader is refered to the paper of Dudley and Norvai\u{s}a \cite{DNo} (see also Gubunelli 
\cite{G}, and Le\'on and Tindel \cite{LT}).

\subsection{Semilinear Volterra integral equations with additive noise}\label{aviq}

Here we consider the Volterra integral equation
\begin{equation}\label{eq1}
X(t)=\xi_t+\frac{1}{\Gamma(\beta)} \int_0^t (t-s)^{\beta-1} A X(s) ds + 
\frac{1}{\Gamma(\alpha)}\int_0^t (t-s)^{\alpha-1}d\theta_s,\quad t\ge 0,
\end{equation}
where the initial condition $\xi=\{\xi_t, t\ge 0\}$ is bounded on compact sets
and measurable, 
$\beta\in(0,1)$, $A\in\mathbb{R}$, $\alpha\in (1,2),$ $\theta=\{\theta_s, s\ge 0\}$ 
is a $\gamma$-H\"older continuous function with $\gamma\in (0,1)$ and $\Gamma$ is the Gamma 
function. The second integral in (\ref{eq1}) is a Young one and it is well-defined if $\alpha-1+
\gamma>1,$ because $s\mapsto (t-s)^{\alpha-1}$ is ($\alpha-1$)-H\"older continuous on $[0,t].$

Now we give two lemmas that we need in the remaining of the paper. The following result 
provides a closed form for the solution to equation (\ref{eq1}) and its proof is in \cite{FLM}.

\begin{lemma}\label{solueq} 
Let $\alpha+\gamma>2$ and $A\in\mathbb{R}.$ Then, the solution to (\ref{eq1}) has the form
\begin{equation}\label{eq2}
X(t)=\xi_t+ A\int_0^t(t-s)^{\beta-1}E_{\beta,\beta}(A(t-s)^\beta)\xi_sds+\int_0^t(t-s)^{\alpha-1}
E_{\beta,\alpha}(A(t-s)^\beta)d\theta_s,\quad t\ge 0.
\end{equation}
Moreover, if
$$\xi_t=\frac{1}{\Gamma(\eta)}\int_0^t(t-s)^{\eta-1}g(s)ds,\quad t\ge 0,$$
where $\eta >0$ and $g\in L^1([0,M))$ for each $M>0,$ then we have
$$X(t)=\int_0^t(t-s)^{\eta-1}E_{\beta,\eta}(A(t-s)^\beta)g(s)ds+\int_0^t(t-s)^{\alpha-1}
E_{\beta,\alpha}(A(t-s)^\beta)d\theta_s,\quad t\ge 0.$$

\end{lemma}

In this work we use comparison methods in order to obtain the stability of some fractional systems. We
can find comparison theorems in the literature for fractional evolution equations (see, e.g. Theorem 4.2 in
\cite{LV}), but, unfortunately this results are not suitable for our purpose. Thus, we give the following lemma, 
that is a version of Theorem 2.2.5 in Pachpatte \cite{Pa} and allows us to prove stability for the semi-linear 
equations that we study. Hence, this result is a fundamental tool in the development of this paper. 

\begin{lemma}\label{pach} Let $k: [0,T]\times \mathbb{R}\to \mathbb{R}$ be a function such that:
\begin{itemize}
\item[i)] $k(\cdot, x)$ is measurable on $[0,T]$ for each $x\in\mathbb{R}.$

\item[ii)]There is a constant $M>0$ such that $|k(s,x)-k(s,y)|\le M|x-y|,$ for any $s\in [0,T]$ and 
$x,y\in\mathbb{R}.$

\item[iii)] $k$ is bounded on bounded sets of $[0,T]\times\mathbb{R}.$

\item[iv)] $k(s,\cdot)$ is non-decreasing for any $s\in [0,T]$.

\end{itemize}
Also, let $B\in\mathbb{R},$ $\beta\in (0,1),$ and $x$ and $y$ two continuous functions on $[0,T]$ such that
\begin{equation}\label{desxy}
x(t)\le y(t)+\int_0^t (t-s)^{\beta-1}E_{\beta,\beta}(B(t-s)^\beta)k(s,x(s))ds,\quad t\in [0,T].
\end{equation}
Then $x\le u$ on $[0,T],$ where $u$ is the solution to the equation
\begin{equation}\label{soluy}
u(t)= y(t)+\int_0^t (t-s)^{\beta-1}E_{\beta,\beta}(B(t-s)^\beta)k(s,u(s))ds,\quad t\in [0,T].
\end{equation}

\end{lemma}
\noindent{\bf Remark}. The assumptions of $k$ yield that equation (\ref{soluy}) has a unique continuous
solution.

\medskip
\noindent{\bf Proof}. Denote by $C([0,T])$ the family of continuous functions on $[0,T].$ 
Let $\mathcal{G}: C([0,T])\to C([0,T])$ be given by

\begin{equation*}
(\mathcal{G}z)(t)=y(t)+\int_0^t (t-s)^{\beta-1}E_{\beta,\beta}(B(t-s)^\beta)k(s,z(s))ds, \quad
t\in [0,T].
\end{equation*}
It is not difficult to see that Hyphoteses i) and iii) imply that $\mathcal{G}$ is well-defined. It means, 
$\mathcal{G}(z)$ is a continuous function for each $z\in C([0,T])$. Remember that we denote $\sup_{t\in[a,b]}|z(t)|$ 
by $||z||_{\infty,[a,b]}$ for any function $z\in C([a,b]).$ Then, from the continuity of $E_{\beta,\beta}$ and 
hypothesis ii), there is a constant $\bar M>0$ such that, for every $z,\tilde z\in C([0,T]),$ we have
\begin{equation*}
\begin{split}
|(\mathcal{G}z)(t)-(\mathcal{G}\tilde z)(t)|&\le \bar M\int_0^t (t-s)^{\beta-1} |k(s,z(s))-k(s,\tilde z(s))|ds
\\ &\le M\bar M\int_0^t (t-s)^{\beta-1} |z(s)-\tilde z(s)|ds
\\ & \le \frac{M\bar M}{\beta}T^\beta ||z-\tilde z ||_{\infty,[0,T]}, \quad\mbox{for}\ t\in [0,T].
\end{split}
\end{equation*}
Similarly, for $\bar T\le T,$ we are able to see that
\[||\mathcal{G}z-\mathcal{G}\tilde z||_{\infty,[0,\bar T]} \le  
\frac{M\bar M}{\beta} ||z-\tilde z ||_{\infty,[0,\bar T]}\bar T^\beta. \]
Consequently, if $\bar T^\beta \frac{M\bar M}{\beta}<1,$ $\mathcal G$ is a contraction on $C([0,\bar T]).$ 
Therefore the sequence $v_{n+1}={\mathcal G} v_n,$ with $v_0=x,$ is such that $v_n(t)\to u(t)$ and $v_n(t)\le 
v_{n+1}(t)$ for $t\in [0,\bar T],$ due to Hypothesis iv),  (\ref{desxy}), and $E_{\beta,\beta}$ being a completely
monotonic function. Thus the result is true if we write $\bar T$ instead
of $T.$

Now, suppose the lemma holds for the interval $[0,n\bar T],$ $n\in\mathbb{N}.$ Then, by (\ref{desxy}) 
we can write
\begin{equation*}
\begin{split}
x(t)&\le y(t)+\int_0^{n\bar T} (t-s)^{\beta-1} E_{\beta,\beta} (B(t-s)^\beta)k(s,x(s))ds +
\int_{n\bar T}^t(t-s)^{\beta-1}E_{\beta,\beta}
(B(t-s)^\beta)k(s,x(s))ds\\
&\le \bar y(t)+\int_{n\bar T}^t (t-s)^{\beta-1}  E_{\beta,\beta} (B(t-s)^\beta)k(s,x(s))ds,\quad
t\in [n\bar{T},(n+1)\bar{T}],
\end{split}
\end{equation*}
where
\[ \bar y(t)= y(t)+\int_0^{n\bar T} (t-s)^{\beta-1} E_{\beta,\beta} (B(t-s)^\beta)k(s,u(s))ds.\]

Finally, defining $\mathcal{G}^{(n)}:C([n\bar T,(n+1)\bar T])\to 
C([n\bar T,(n+1)\bar T])$ by
\[ ({\mathcal G^{(n)}} z)(t)=\bar y(t)+\int_{n\bar T}^{t}(t-s)^{\beta-1}
E_{\beta,\beta}(B(t-s)^{\beta})k(s,z(s))ds, \]
and using the fact that equation (\ref{soluy}) has a unique solution due to Hypothesis ii),
we can proceed as in the first part of this proof to see that $x\le u$ on $[0,(n+1)\bar{T}]$. Thus, 
the result follows using induction on $n$.

\hfill $\Box$

\section{A class of a nonlinear fractional-order systems}\label{clasemi1}

In this section we establish two sufficient conditions for the stability of
a deterministic semilinear Volterra integral equation. Thus, we improve the results in \cite{FLM}
for this kind of systems when the noise is null (i.e., $Z$ in (\ref{maineq}) is equal to zero). 

\subsection{A constant as initial condition}\label{concons}

This part is devoted to refine Theorem 1 of \cite{WWL} in the one-dimensional case.
Toward this end, in this section, we suppose that the initial condition is a constant. That is,
we first consider the fractional equation
\begin{equation}\label{equ3} 
X(t)=x_0+\frac{1}{\Gamma(\beta)} \int_0^t (t-s)^{\beta-1} A X(s) ds
 + \frac{1}{\Gamma(\beta)} \int_0^t (t-s)^{\beta-1} h(X(s)) ds,\quad t\ge 0,
\end{equation}
with $x_0\in \mathbb{R}$, $\beta\in(0,1)$, $A<0$ and $h:\mathbb{R}\to
\mathbb{R}$ a measurable function.

In the remaining of this paper we deal with the following hypotheses.
\begin{itemize}
\item[(H1)]  There is a constant $C>0$ such that $A+C<0$ and $|h(x)|\le C|x|$, for all $x\in \mathbb{R}.$
\item[(H2)]  There are $\delta_0>0$ and $C>0$ such that $A+C<0$ and $|h(x)|\le C|x|,$ for $|x|<\delta_0.$ 
\end{itemize}

Now, we consider several definitions of stability.

\begin{definition}\label{glar}
Any solution $X$ to equation (\ref{equ3}) is said to be:
\begin{itemize}

\item[i)] {\rm globally stable in the large} if $X(t)$ goes to zero as $t$ tends to infinity, 
for all $x_0\in\mathbb{R}$ . 

\item[ii)] {\rm Mittag-Leffler stable} if there is $\delta>0$ such that $|x_0|<\delta$
 implies
$$|X(t)|\le \left[ m(x_0)E_{\beta,1}(Bt^\beta) \right]^b,\quad t\ge 0,$$
where $\beta\in (0,1),$ $B<0,$ $b>0$ and $m$ is a positive and locally Lipschitz function with $m(0)=0.$

\item[iii)] {\rm stable} if for $\varepsilon>0,$ 
there is $\delta>0$ such that $|x_0|<\delta$ implies $|X(t)|<\varepsilon,$ \ for all $t\ge0.$ 

\item[iv)] {\rm stable in the large} if there is $\delta>0$ such that $|x_0|<\delta$ 
implies $\lim_{t\to\infty}X(t)=0.$

\item[v)]  {\rm asymptotically stable} if it is stable and stable in the large.

\end{itemize}
\end{definition}

\begin{remark}\label{remsolu} Observe that, under the assumptions that $h$ is continuous and satisfies 
(H1), equation (\ref{equ3}) has at least one solution on $[0,\infty)$
because of \cite{LV} (Theorems 3.1 and 4.2). Indeed, in \cite{LV} (Theorem 4.2) we can consider
\[ g(t,x) = \left\{ \begin{array}{ll}
         0 & \mbox{if $x \le 0$},\\
        (|A|+C)x & \mbox{if $x>0$}.\end{array} \right. \]       
Similarly, for a continuous function $h$ satisfiying (H2), we introduce the function

\[ \varphi(x) = \left\{ \begin{array}{ll}
         x & \mbox{if $|x| \le \delta_0 /2$},\\
        \delta_0 & \mbox{if $x>\delta_0 /2$},\\
   -\delta_0 & \mbox{if $x<-\delta_0 /2$}      .\end{array} \right. \]     
Then, using \cite{LV} (Theorems 3.1 and 4.2) again, the equation
\begin{equation}\label{xaux}
X(t)=x_0+\frac{1}{\Gamma(\beta)}\int_0^t(t-s)^{\beta-1}(AX(s)+h(\varphi(X(s)))ds
\end{equation}
has at least one solution defined on $[0,\infty)$ due to $|Ax+h(\varphi(x))|\le |Ax|+C|\varphi(x)|
\le (|A|+C)|x|.$ Hence equation (\ref{equ3}) has at least one continuous solution on $[0,\infty)$ if (\ref{xaux}) 
is stable and $x_0$ is small enough because, in this case, the solution of (\ref{xaux}) is also a solution of 
equation (\ref{equ3}) and $h\circ \varphi$ is bounded. So, without loss of generality we can assume that 
(\ref{equ3}) has at least one continuous solution because one of the main purposes of the paper is to deal with 
the stability of (\ref{maineq}).
\end{remark}

We need the following lemma to prove some of our results. The main idea of its proof is in the paper of Mart\'inez-Mart\'inez 
et al. \cite{MM-L-FA}. Here we give an sketch of the proof for the convenience of the reader. 

\begin{lemma}\label{rafa}
Let $h$ be as in (H2) (resp. (H1)). Then, for $0<x_0<\delta_0$ (resp. $x_0>0$), any continuous solution $X$
of (\ref{equ3}) satisfies $X(t)>0$ for all $t\ge 0.$
\end{lemma}
\noindent{\bf Proof} (An idea). Let (H2) (resp. (H1)) be true and $x_0\in (0,\delta_0)$ (resp. $x_0>0$). 
Then the continuity of $X$ implies that there is $\tau>0$ such that $X(t)\in (0,\delta_0)$ 
(resp. $X(t)>0$) for all $t\in [0,\tau].$ Consequently
\begin{equation}\label{igzero1}
0<X(t)\le x_0+\frac{1}{\Gamma(\beta)}\int_0^t (t-s)^{\beta-1}[A+C]X(s)ds<x_0,\quad t\in [0,\tau].
\end{equation}
In other words, we have proved that $X(t)$ is less than $x_0$ if $X>0$ on $[0,t].$
Now suppose that 
$\tau_0=\inf \{ t>0: X(t)=0 \}$ is finite. Hence, from (\ref{equ3}) we deduce
\begin{equation*}
x_0=-\frac{1}{\Gamma(\beta)}\int_0^{\tau_0}(\tau_0-s)^{\beta-1}[AX(s)+h(X(s))]ds.
\end{equation*}
Therefore, using (\ref{equ3}) again, we have
\begin{equation}\label{igu0}
X(t)=\frac{1}{\Gamma(\beta)}\left(\int_0^{t}(t-s)^{\beta-1}[AX(s)+h(X(s))]ds
-\int_0^{\tau_0}(\tau_0-s)^{\beta-1}[AX(s)+h(X(s))]ds\right),
\quad t\le\tau_0.
\end{equation}

Since $AX(s)+h(X(s))\le [A+C]X(s)\le 0$ on $[0,\tau_0]$ due to the hipothesis of this result, we 
are able to write
\begin{equation}\label{descot}
X(t)\le\frac{1}{\Gamma(\beta)}(|A|+C)\int_t^{\tau_0}(\tau_0-s)^{\beta-1}X(s)ds,\quad t\in [0,\tau_0].
\end{equation}
So, as $(\tau_0-s)^{-1}\le (t-s)^{-1},$ $X(\tau_0)=0$ and the continuity of the solution of equation (\ref{equ3}),
iterating inequality (\ref{descot}) we can find a positive constant $\tilde{C}$ such that
$$X(t)\le \tilde{C} 
\left( \frac{(|A|+C)(\tau_0-t)^\beta}{\Gamma(\beta+1)}\right)^{n-1}.$$ 
Thus, taking $t\in (0, \tau_0)$ such that $\frac{(|A|+C)(\tau_0-t)^\beta}{\Gamma(\beta+1)}<1$ we deduce 
$X(t)=0,$ which is a contradiction.

\hfill $\Box$

\medskip

\noindent{\bf Remark}.
As it was pointed out in \cite{MM-L-FA}, if the initial condition in equation (\ref{equ3}) is a non-decreasing, continuous and non-negative
function instead of a constant, we can repeat the procedure in this proof in order to obtain the same result. 
Indeed, suppose that $0<\xi_t<\delta_0$ (resp. $\xi_t> 0$) for all $t\ge 0$,
first of all (\ref{igzero1}) becomes
$$ 0<X(t)\le \xi_t+\frac{1}{\Gamma(\beta)}\int_0^t (t-s)^{\beta-1}[A+C]X(s)ds<\xi_t.$$
Secondly, instead of equallity (\ref{igu0}) we have
\begin{equation*}
\begin{split}
X(t)&=\xi_t-\xi_{\tau_0}+\frac{1}{\Gamma(\beta)}\left(\int_0^{t}(t-s)^{\beta-1}[AX(s)+h(X(s))]ds
-\int_0^{\tau_0}(\tau_0-s)^{\beta-1}[AX(s)+h(X(s))]ds\right)\\ &\le 
\frac{1}{\Gamma(\beta)}\left(\int_0^{t}(t-s)^{\beta-1}[AX(s)+h(X(s))]ds
-\int_0^{\tau_0}(\tau_0-s)^{\beta-1}[AX(s)+h(X(s))]ds\right),
\quad t\le\tau_0,
\end{split}
\end{equation*}
due to $\xi$ being non-decreasing. Hence, it is not difficult to see that (\ref{descot}) is still satisfied. 

\medskip

An immediate consequence of the first part of the proof of Lemma \ref{rafa} is the following.

\begin{corollary}
Assume either (H2), or (H1) is satisfied. Then, any continuous solution to equation (\ref{equ3}) is stable. 
\end{corollary}
\noindent{\bf Proof}. If $x_0>0,$ the result follows from (\ref{igzero1}).

For $x_0<0$ and $X$ a solution of (\ref{equ3}), we have $-X$ is a solution of  
\begin{displaymath}
Y(t)=-x_0+\frac{1}{\Gamma(\beta)} \int_0^t (t-s)^{\beta-1} [A  Y(s) ds + \hat{h}(Y(s))] 
ds,\quad t\ge 0,\end{displaymath}
with $\hat{h}(x)=-h(-x).$
\hfill $\Box$

Now we establish the main result of this subsection.

\begin{proposition}\label{P311} Let $h$ be a function satisfiying (H2) (resp. (H1)). Then, any continuous solution 
of equation (\ref{equ3}) is Mittag-Leffler stable and therefore is also asymptotically stable
(resp. globally stable in the large).
\end{proposition}

\noindent {\bf Proof}. Let (H2) (resp. (H1)) be satisfied and $0<x_0<\delta_0$ 
(resp. $x_0>0$). Then $0< X(t)<\delta_0$ (resp. $X(t)>0$) by Lemma \ref{rafa} and its proof (see (\ref{igzero1})). 

On the othe hand, consider the solution $Z$ of the following linear fractional equation
\begin{displaymath}
Z(t)=2x_0+\frac{1}{\Gamma(\beta)} \int_0^t (t-s)^{\beta-1} [A+C] Z(s) ds,\quad t\ge 0.\end{displaymath}
Then by the continuity of the solutions $X$ and $Z$, there exists $\tau>0$ such that, for all $t\in (0,\tau)$, 
we have $0<X(t)<Z(t)$. 
If this inequality is satisfied for any $t>0$, we can ensure that 
$X$ is asymptotically stable (resp. and globally stable in the large), and that this solution also is Mittag-Leffler stable because 
the solution to $Z$ of last equation is given by (see \cite{JJM} or Lemma \ref{solueq})
$$Z(t)=2x_0E_{\beta,1}([A+C]t^\beta),\quad t\ge 0.$$

We now suppose that 
there exists  $t_0>0$ such that $X(t_0)=Z(t_0)$ and $X(t)<Z(t),$ for $t<t_0.$ Set $Y=X-Z$, then
\begin{displaymath}
Y(t)=-x_0+\frac{1}{\Gamma(\beta)} \int_0^t (t-s)^{\beta-1} A Y(s) ds + 
\frac{1}{\Gamma(\beta)} \int_0^t (t-s)^{\beta-1} [h(X(s))-C Z(s)] ds,\quad t\ge 0.\end{displaymath}
From (\ref{eq2}) (see also \cite{JJM}) we observe that $Y$ also satisfies the equality
\begin{displaymath}
Y(t)=-x_0 E_{\beta,1}(At^\beta) + \int_0^t (t-s)^{\beta-1}
E_{\beta,\beta}(A(t-s)^\beta) [h(X(s))-C Z(s)] ds,\quad t\ge 0.\end{displaymath}
For $s\in (0,t_0)$, we have $|h(X(s))|\le CX(s)< C Z(s).$ Thus $h(X(s))-C Z(s)<0$. Consequently, by
the completely monotonic property of $E_{\beta,\beta}$ we have $Y(t_0)<0,$ and this is a 
contradiction because it is supposed that $Y(t_0)=0$. Now we can conclude that $X$ is Mittag-Leffler stable. 

Finally we consider the case that $-\delta_0<x_0<0$ (resp. $x_0<0$). Note that $\hat X=-X$ is such that
\begin{displaymath}
\hat X(t)=-x_0+\frac{1}{\Gamma(\beta)} \int_0^t (t-s)^{\beta-1} A \hat X(s) ds + 
\frac{1}{\Gamma(\beta)} \int_0^t (t-s)^{\beta-1} \tilde h(\hat X(s)) ds,\quad t\ge 0,\end{displaymath}
with $\tilde h(x)=-h(-x)$. Hence, by the first part of this proof and the fact that $\tilde{h}$ satisfies
(H2) (resp. (H1)), we have that the proof is complete.

\hfill $\Box$

\medskip

\noindent {\bf Remark}. Let $X$ be a solution to equation (\ref{equ3}). 
Wen et al. \cite{WWL} (Theorem 1) have  proved that the solution to equation 
(\ref{equ3}) is stable if $\lim_{|x|\to 0}\frac{|h(x)|}{|x|}\to 0.$ 
Also, Zhang and Li \cite{ZC} have used a result similar to Lemma \ref{solueq} to prove that $X$
is asymptotically stable for the case that $\lim_{x\to 0}\frac{|h(x)|}{|x|}=0, \beta\in (1,2)$
and $\beta+\frac{1}{|A|}<2.$ 
Proposition \ref{P311} establishes that $X$ 
is asymptotically stable under a weaker condition. Namely (H2). This is possible because we use a 
comparison type result and the fact that this solution does not change sign. 

\subsection{A function as initial condition}\label{ficon}

Here we treat the case that the initial condition is a function satisfying some suitable conditions. 

Consider the following deterministic Volterra integral equation
\begin{equation}\label{eqspa}
X(t)=\xi_t+\frac{1}{\Gamma(\beta)} \int_0^t (t-s)^{\beta-1} A X(s) ds + 
\frac{1}{\Gamma(\beta)} \int_0^t (t-s)^{\beta-1} h(X(s)) ds,\quad
t\ge 0.\end{equation}
Here $\beta\in(0,1)$, $A<0$, and $h:\mathbb{R}\longmapsto \mathbb{R}$ and $\xi:\mathbb{R}^+
\longmapsto \mathbb{R}$ are two measurable functions.

Concerning the existence of a continuous solution of equation (\ref{eqspa}) we remark the following. 
For a continuous function $h$ as in (H1) and $\xi$ continuous,
 we can consider the equation 
$$Z(t)=\frac{1}{\Gamma(\beta)}\int_0^t (t-s)^{\beta-1}f(s,Z(s))ds,$$
where $f(s,x)=A(x+\xi_s)+h(x+\xi_s),$ which has a solution $Z$ due to Theorem 4.2 in \cite{LV} 
(with $g(s,x)=(|A|+C)(x+|\xi_s|)$) and Lemma \ref{solueq}. Therefore $Z+\xi$ is a 
solution of (\ref{eqspa}). Similarly if $\xi$ is ``small enough" and $h$ is either a continuous Lipschitz
function on a neighbourhood of zero, or as in (H2), then we can proceed as in Remark \ref{remsolu} 
to see that (\ref{eqspa}) has at least one solution in this case. Therefore, as in Remark \ref{remsolu},
we can assume that (\ref{eqspa}) has at least one continuous solution.  


\bigskip

On the other hand, in this paper we analyze several stability criteria for different classes $\mathcal{E}$ of initial 
conditions. Sometimes $\mathcal{E}$ is a subset of a normed linear space $\mathcal{X}$ of 
continuous functions endowed with the norm $||\cdot||_{\mathcal{X}}.$ In other
words we consider normed linear spaces $(\mathcal{X},||\cdot||_{\mathcal{X}}).$ Mainly, 
in the remaining of this paper, we deal with the following classes of initial conditions.
 
\begin{definition}\label{coninic} We have the following assumptions on $\xi$:

\begin{enumerate}
\item[1.] If the initial condition $\xi$ is continuous on $[0,\infty)$ and there is $\xi_\infty\in\mathbb{R}$ 
such that, given $\varepsilon>0$, there exists $t_0>0$ such that $|\xi_s-\xi_\infty|\le \varepsilon$ 
for any $s\ge t_0,$ we say that $\xi$ belongs to the family $\mathcal{E}^1.$

\item[2.] $\mathcal{E}^2$ is the set of all functions $\xi$ of class ${C}^1(\mathbb{R}_+)$
(i.e. $\xi$ has a continuous derivative on $\mathbb{R}_+$) such that
\begin{equation*}\lim_{t \mapsto \infty} |\xi_t|/t^\beta=0 \quad 
\mbox{and} \quad  |\xi_t'| \le \frac{\tilde C}{t^{1-\upsilon}}, \quad \textrm{for some}\ 
\upsilon \in (0,\beta) \ \textrm{and}\ \tilde C\in\mathbb{R}. \end{equation*}

\item[3.] $\mathcal{E}^3$ is the space of continuous functions of the form
\begin{equation}\label{xig}
\xi_t=\frac{1}{\Gamma(\eta)}\int_0^t (t-s)^{\eta-1} g(s) ds,
\end{equation}
with $g \in L^1([0,\infty)) \cap L^p([0,\infty)),$ $\eta \in (0,\beta+1)$ and
$p>\frac{1}{\eta}\vee 1$. 
\end{enumerate}
\end{definition}

The stability concepts that we develop in this section are the following. 
\begin{definition}\label{gsle} Let $\mathcal{E}\subset\mathcal{X}.$ A solution $X$ of (\ref{eqspa})
is said to be:
\begin{itemize}
\item[i)]  {\rm globally stable in the large for 
the class} $\mathcal{E}$ (or {\rm globally $\mathcal{E}$-stable} in the large) if $X(t)$ tends to zero as 
$t\to\infty,$ \ for every $\xi\in\mathcal{E}.$

\item[ii)] $\mathcal{E}$-{\rm stable} if for $\varepsilon>0,$ there is $\delta>0$ such that $||X||_{\infty,[0,\infty)}
<\varepsilon$ for every $\xi\in\mathcal{E}$ satisfiying $||\xi||_{\mathcal{X}}<\delta.$  

\item[iii)] asymptotically $\mathcal{E}$-stable if it is $\mathcal{E}$-stable and there is $\delta>0$ such that
$\lim_{t\to\infty}X(t)=0$ for any $\xi\in\mathcal{E}$ such that $||	\xi||_{\mathcal{X}}<\delta.$
\end{itemize}

\end{definition}

In the following auxiliary result, $\mathcal{E}^4$ is the family of functions $\xi$ having the
form (\ref{xig}) with $\eta=\beta$ and $g$ is a continuous function such that $\lim_{t\to\infty}
g(t)=0.$ In this case, the involved norm is $||\xi||_{\mathcal{X}}=||g||_{\infty,[0,\infty)}.$

\begin{lemma}\label{remarkag}  Let $B<0$ and $\xi\in\mathcal{E}^4.$
Then the solution to the equation
$$Y(t)=\xi_t+\frac{1}{\Gamma(\beta)}\int_0^t(t-s)^{\beta-1}BY(s)ds,\quad  t\ge 0,$$
is $\mathcal{E}^4$-stable and globally $\mathcal{E}^4$-stable in the large.
\end{lemma}
\noindent{\bf Proof}. We observe that, by Lemma \ref{solueq}, we have
\begin{equation*}
\begin{split}
Y(t)&=\int_0^t (t-s)^{\beta-1}E_{\beta,\beta}(B(t-s)^\beta)g(s)ds= 
\int_0^t s^{\beta-1}E_{\beta,\beta}(Bs^\beta)g(t-s)ds,\quad t\ge 0.
\end{split}
\end{equation*}
So, the completely monotone property of $E_{\beta,\beta},$ (\ref{estiml}) and (\ref{intmitlef}) lead us 
to establish 

\begin{equation*}
\begin{split}
|Y(t)|&\le \left(\sup_{s\ge 0} |g(s)|\right)\int_0^ts^{\beta-1}E_{\beta,\beta}(Bs^\beta)ds=
\left(\sup_{s\ge 0} |g(s)|\right)t^{\beta}E_{\beta,\beta+1}(Bt^\beta)\\
&\le \frac{C_{\beta,\beta+1}}{|B|} ||g||_{\infty,[0,\infty)}.
\end{split}
\end{equation*}
Thus, $Y$ is $\mathcal{E}^4$-stable.

Also, by (\ref{intmitlef}) we are able to write
\begin{equation*}
\begin{split}
Y(t)&=\int_0^t (t-s)^{\beta-1}E_{\beta,\beta}(B(t-s)^\beta)g(s)ds\\
&= g(t)t^\beta E_{\beta,\beta+1}(Bt^\beta)+\int_0^t (t-s)^{\beta-1}E_{\beta,\beta}(B(t-s)^\beta)
[g(s)-g(t)]ds,\quad t\ge 0.
\end{split}
\end{equation*}
Therefore, using (\ref{estiml}) and the proof of Proposition 3.3.1 in \cite{FLM} again, together with the facts that 
$B<0$ and $g$ is a continuous function such that $\lim_{t\to\infty}g(t)=0,$ we obtain $Y(t)\to 0$ as $t\to\infty.$

\hfill $\Box$

Now we give a general result. 

\begin{theorem}\label{teolicon} Let (H2) (resp. (H1)) be true, and $\mathcal{E}$ a family of continuous
functions of a normed linear space $\mathcal{X}$ such that the solution of the equation
\begin{equation}\label{volt1}
Y(t)=\xi_t+\frac{1}{\Gamma(\beta)}\int_0^t(t-s)^{\beta-1}AY(s)ds,\quad t\ge 0,
\end{equation}
is asymptotically $\mathcal{E}$-stable (resp. globally 
$\mathcal{E}$-stable in the large). Then any continuous solution of equation (\ref{eqspa}) is also
asymptotically $\mathcal{E}$-stable (resp. globally $\mathcal{E}$-stable in the large).
\end{theorem}
\noindent{\bf Proof}. Suppose that (H1) (resp. (H2)) is true. Let $X$ be a continuous solution to equation
(\ref{eqspa}). Take $Z=X-Y,$ then we have 
\begin{equation*}
Z(t)=\frac{1}{\Gamma(\beta)}\int_0^t (t-s)^{\beta-1}[AZ(s)+h(X(s))]ds,\quad t\ge 0.
\end{equation*}
Thus, Lemma \ref{solueq} allows us to write 

\begin{equation*}
Z(t)=\int_0^t (t-s)^{\beta-1} E_{\beta,\beta}(A(t-s)^\beta)h(X(s))ds,\quad t\ge 0.
\end{equation*}
Hence, for $\xi\in\mathcal{E}$ we have (resp. for $\xi\in\mathcal{E}$ such that $||Y||_{\infty,[0,\infty)}<
\delta_0,$ which gives $|\xi_0|=|Y(0)|<\delta_0,$ the continuity of $X$ implies that there is $t_0>0$
such that $||X||_{\infty,[0,t_0)}<\delta_0$ and)
\begin{equation*}
\begin{split}
|Z(t)|&\le C\int_0^t (t-s)^{\beta-1} E_{\beta,\beta}(A(t-s)^\beta)|X(s)|ds\\
 &\le C\int_0^t (t-s)^{\beta-1} E_{\beta,\beta}(A(t-s)^\beta)|Z(s)|ds+
 C\int_0^t (t-s)^{\beta-1} E_{\beta,\beta}(A(t-s)^\beta)|Y(s)|ds,\ t\ge 0\ (\mbox{resp.} \ t\le t_0),
\end{split}
\end{equation*}
where we make use of the completely monotonic property of $E_{\beta,\beta}.$
Invoking Lemma \ref{pach} and the uniqueness of the solutions for the involved equations, 
$|Z(t)|\le u(t)$ for all $t\ge 0$ (resp. $t\le t_0$), where $u$ is the solution to
\begin{equation*}
u(t)=\frac{C}{\Gamma(\beta)}\int_0^t (t-s)^{\beta-1}|Y(s)|ds+
\frac{1}{\Gamma(\beta)}\int_0^t (t-s)^{\beta-1} [A+C]u(s)ds,\quad t\ge 0.
\end{equation*}
Finally observe that $|X(t)|\le u(t)+|Y(t)|$ for $t\ge 0$ (resp. for $t\le t_0$ such that $||X||_{\infty, [0,t_0)}
<\delta_0$). Thus Lemma \ref{remarkag} implies that $u$ is globally $\mathcal{E}^4$-stable in the large 
(resp. $u$ is $\mathcal{E}^4$-stable and globally $\mathcal{E}^4$-stable in the large), wich gives that
the proof is complete. 

\hfill $\Box$

\noindent{\bf Remark}. For each $i\in \{1,\ldots,n\}$ let $\mathcal{X}^{i}$ be a normed linear space
of functions. Note that if $\xi=\sum_{i=1}^n\xi^{(i)},$ where $\xi^{(i)}\in\hat{\mathcal{E}}^{i}
\subset\mathcal{X}^{i}$ and (\ref{volt1}) is $\hat{\mathcal{E}}^{(i)}$-stable for each $i\in \{1,\ldots,n\}.$
Then, (\ref{volt1}) is also $\mathcal{E}$-stable, where $\mathcal{E}$ is the family of functions of the form
$\sum_{i=1}^n\xi^{(i)}$ and the involved seminorm is $||\xi||_{\mathcal{X}}=\sum_{i=1}^n||\xi^{(i)}
||_{\mathcal{X}^{i}}.$ Indeed, by Lemma \ref{solueq} we have that the solution $Y$ is given by 
$$Y(t)=\sum_{i=1}^nY^{(i)}(t)=\sum_{i=1}^n\left(\xi_t^{(i)}+A\int_0^t(t-s)^{\beta-1}E_{\beta,\beta}
 (A(t-s)^\beta)\xi_s^{(i)}ds  \right),
\quad t\ge 0,$$
where, for each $i\in \{1,\ldots, n\},$ $Y^{(i)}$ is the unique solution to the linear equation
\[ Y^{(i)}(t)=\xi_t^{(i)}+\frac{1}{\Gamma(\beta)}\int_0^t(t-s)^{\beta-1}AY^{(i)}(s)ds,\quad t\ge 0. \]

In the following result we see that the family $\mathcal{E}:=\{ \xi\in C([0,\infty)): \xi=\sum_{i=1}^3\xi^{(i)},\xi^{(i)}\in
\mathcal{E}^i\}$ is an example of a family of functions for which the assumptions of Theorem \ref{teolicon} is
satisfied. Here, $||\cdot||_{\mathcal{X}^{1}}=
||\cdot||_{\infty,[0,\infty)},$ $||\xi^{(2)}||_{\mathcal{X}^{2}}=||\xi_{\cdot}^{(2)}E_{\beta,1}(A\cdot^{\beta})||_{\infty,[0,\infty)}$
$+||\cdot^{1-\upsilon}{\xi_{\cdot}^{(2)}}'||_{\infty,[0,\infty)}$  and $||\xi^{(3)}||_{\mathcal{X}^{3}}=||g||_{L^1([0,\infty))}$
+$||g||_{L^p([0,\infty))},$ where  $\cdot^{1-\upsilon}{\xi_{\cdot}^{(2)}}'$ denotes 
$s\mapsto s^{1-\upsilon}{\xi^{(2)}_s}'$ and $\xi^{(3)}$ is given by the right-hand side of (\ref{xig}). Thus, in this case 
$||\xi||_{\mathcal{X}}=\sum_{i=1}^3 ||\xi^{(i)}||_{\mathcal{X}^{i}}.$

\begin{proposition}\label{ejemcomb} Let $A<0$ and $\beta\in (0,1).$ Then any solution to 
(\ref{volt1}) is $\mathcal{E}$-stable and $\mathcal{E}$-stable in the large.
\end{proposition}
\noindent {\bf Proof}. By previous remark we only need that equation
(\ref{volt1}) is $\mathcal{E}^{i}$-stable and $\mathcal{E}^{i}$-stable in the large, for $i=1,2,3$. To prove this, let
$Y$ be the solution to equation (\ref{volt1}). The global $\mathcal{E}^{i}$-stability in the large has already been considered
in \cite{FLM} (Theorem 3.3). Now we divide the proof in three steps.

\medskip

\noindent{\bf Step 1}. Here we consider the case $i=1.$ Then Lemma \ref{solueq} and (\ref{intmitlef}) give that, 
for $t\ge 0,$
\begin{equation*}
\begin{split}
|Y(t)| &\le |\xi^{(1)}_t|+|A|\int_0^t (t-s)^{\beta-1}E_{\beta,\beta} (A(t-s)^\beta)|\xi_s^{(1)}|ds\\
&\le ||\xi^{(1)}||_{\infty,[0,\infty)}\left( 1+|A|\int_0^t (t-s)^{\beta-1}E_{\beta,\beta} (A(t-s)^\beta)ds\right)\\
&= ||\xi^{(1)}||_{\infty,[0,\infty)}\left( 1+|A|t^\beta E_{\beta,\beta+1}(At^\beta)\right)\\
&\le ||\xi^{(1)}||_{\infty,[0,\infty)}(1+C_{\beta,\beta+1}),
\end{split}
\end{equation*}
which implies that the solution of (\ref{volt1}) is $\xi^{(1)}$-stable.

\medskip

\noindent{\bf Step 2}. For $i=2,$ we get
\begin{equation*}
\begin{split}
|Y(t)|&\le |\xi^{(2)}_tE_{\beta,1}(At^\beta)|+\Big{|}A\int_0^t (t-s)^{\beta-1}E_{\beta,\beta}(A(t-s)^\beta)
(\xi^{(2)}_s-\xi^{(2)}_t)ds \Big{|} \\
&\le ||\xi^{(2)}||_{\mathcal{X}^{(2)}}\left( 1+|A|\int_0^t (t-s)^{\beta}E_{\beta,\beta} (A(t-s)^\beta)
s^{\upsilon-1}ds\right),\quad t\ge 0.
\end{split}
\end{equation*}
Consequently, \cite{FLM} (proof of Theorem 3.2.2) yields 
\begin{equation*}
\begin{split}
|Y(t)| &\le ||\xi^{(2)}||_{\mathcal{X}^{(2)}}\left( 1+t^\upsilon\Gamma(\upsilon)[\upsilon E_{\beta,\upsilon+1}
(At^\beta)-E_{\beta,\upsilon}(At^\beta)] \right)\\
&\le C||\xi^{(2)}||_{\mathcal{X}^{(2)}},\quad t\ge 0,
\end{split}
\end{equation*}
where $C>0$ is a constant and we have utilized that $\upsilon<\beta.$

\medskip

\noindent{\bf Step 3}. Finally we consider the case $i=3.$ In this scenario, 
from Lemma \ref{solueq}, we obtain
\begin{equation*}
\begin{split}
|Y(t)|&= \Big{|} \int_0^t (t-s)^{\eta-1}E_{\beta,\eta}(A(t-s)^\beta)g(s)ds \Big{|}\\
&= \Big{|} \int_0^t s^{\eta-1}E_{\beta,\eta}(As^\beta)g(t-s)ds \Big{|}\\
&\le \Big{|} \int_0^{t\wedge 1} s^{\eta-1}E_{\beta,\eta}(As^\beta)g(t-s)ds \Big{|}
+\Big{|} \int_{t\wedge 1}^t s^{\eta-1}E_{\beta,\eta}(As^\beta)g(t-s)ds \Big{|}\\
&= I_1^{(3)}(t)+I_2^{(3)}(t),\quad t\ge 0.
\end{split}
\end{equation*}
For $I_1^{(3)}$ we can apply H\"older inequality to write, for $q^{-1}=1-p^{-1}$
and $C>0$
\begin{equation*}
\begin{split}
I_1^{(3)}(t)&\le C_{\beta,\eta}\left[ \int_0^1s^{q(\eta-1)}ds \right]^{1/q}
\left[  \int_0^{t\wedge 1}|g(t-s)|^p ds\right]^{1/p}\\
&\le C||g||_{L^p([0,\infty))},\quad t\ge 0,
\end{split}
\end{equation*}
and for $I_2^{(3)}$ we use the fact that $\eta-1-\beta<0.$ Thus
$$|I_2^{(3)}(t)|\le \frac{C_{\beta,\eta}}{|A|}||g||_{L^1([0,\infty))}.$$

\hfill $\Box$

\medskip

\noindent {\bf Remark}. Observe that $\mathcal{E}^1$ contains the bounded variation functions 
on compact sets of $\mathbb{R}_+$  of the form $\xi=\xi^{(1)}-\xi^{(2)},$ where $\xi^{(1)}$ and $\xi^{(2)}$ are two
non-decreasing and bounded functions on $\mathbb{R}_+.$ 

\medskip

The following result is an immediate consequence of Theorem \ref{teolicon} and Proposition
\ref{ejemcomb}.

\begin{theorem}\label{tecor} Suppose that (H2) (resp. (H1)) holds. Let $\xi$ be as in Proposition
\ref{ejemcomb}. Then, any continuous solution to (\ref{eqspa}) is asymptotically $\mathcal{E}$-stable
(resp. globally $\mathcal{E}$-stable in the large).

\end{theorem}

\bigskip

\section{Semilinear integral equations with additive noise}\label{adino1}

In  this section we consider the equation 
\begin{equation}\label{eqadn}
X(t)=\xi_t+\frac{1}{\Gamma(\beta)}\int_0^t (t-s)^{\beta-1} [AX(s)+h(X(s))]ds+
\frac{1}{\Gamma(\alpha)}\int_0^t (t-s)^{\alpha-1}f(s)d\theta_s, \quad t\ge 0.
\end{equation}
Here $\xi,\beta,A$ and $h$ are as in equation (\ref{eqspa}). Henceforth we assume that
$\alpha\in (1,2),$ $\theta=\{ \theta_s, s\ge 0\}$ is a $\gamma$-H\"older continuous
function with $\gamma\in (0,1)$ such that $\theta_0=0$ and $\gamma+\alpha>2,$ and
$f$ is a $\tau$-H\"older continuous function in $C^1(\mathbb{R}_+),$ with $\tau+\gamma>1.$
Note that, in this case, the Young integral in the right-hand side of (\ref{eqadn}) is equal to
$\frac{1}{\Gamma(\alpha)}\int_0^t (t-s)^{\alpha-1}d\tilde{\theta_s},$ where
$\tilde{\theta_s}=\int_0^sf(r)d\theta_r$ due to \cite{FLM} (Lemma 2.4). Thus, Lemma \ref{solueq}
is still true for (\ref{eqadn}) and \cite{FLM} (Lemma 2.7) implies 
$$\frac{1}{\Gamma(\alpha)}\int_0^t (t-s)^{\alpha-1}f(s)d\theta_s=\frac{\alpha-1}{\Gamma(\alpha)}
\int_0^t (t-s)^{\alpha-2}\tilde{\theta}_sds.$$
Hence, the existence of a continuous solution to (\ref{eqadn}) can be considered as in Section \ref{ficon}.

\begin{definition}\label{estdefru}
Let $\mathcal{E}\subset\mathcal{X}$ be a family of continuous functions. We say that a solution 
$X$ of (\ref{eqadn}) is
\begin{enumerate}
\item[i)] ($\mathcal{E},p$)-stable if for $\varepsilon>0,$ there is $\delta>0$ such that
$||X||_{\infty,[0,\infty)}<\varepsilon$ for any $(\xi,f,\theta)$ such that
\begin{equation}\label{a1}
||\xi||_{\mathcal{X}}+||f\theta||_{L^1([0,\infty))}+||f\theta||_{L^p([0,\infty))}+
||\dot f\theta||_{L^1([0,\infty))}<\delta.
\end{equation}


\item[ii)] {\rm asymptotically ($\mathcal{E},p$)-stable } if it is ($\mathcal{E},p$)-stable and there is
$\delta>0$ such that $\lim_{t\to\infty} X(t)=0$ for any ($\xi,f,\theta$) satisfiying (\ref{a1}).
\end{enumerate}
\end{definition}

An extension of Theorem \ref{teolicon} is the following. 

\begin{theorem}\label{sim2} 
Let (H2) (resp. (H1)) be satisfied and $\mathcal{E}$ a class of continuous functions such that
the solution of the equation 
\begin{equation}\label{n42}
Y(t)=\xi_t+\frac{1}{\Gamma(\beta)}\int_0^t(t-s)^{\beta-1}AY(s)ds+
\frac{1}{\Gamma(\alpha)}\int_0^t (t-s)^{\alpha-1}f(s)d\theta_s,\quad t\ge 0,
\end{equation}
is asymptotically ($\mathcal{E},p$)-stable (resp. globally $\mathcal{E}$-stable in the large).
Then, any continuous solution of (\ref{eqadn}) is also 
asymptotically ($\mathcal{E},p$)-stable (resp. globally $\mathcal{E}$-stable in the large).
\end{theorem}
\noindent{\bf Proof}. Observe $X(0)=\xi_0.$ Consequently the proof is similar to that of 
Theorem \ref{teolicon}.

\hfill $\Box$

\medskip

Now we state a consequence of Theorem \ref{sim2}.
\begin{theorem} \label{sim3} Assume  (H2) (resp. (H1)) holds. Let $\xi$ be as in Proposition
\ref{ejemcomb}, $f\in C^1((0,\infty))$ such that $\dot f\theta\in L^1([0,\infty))$ and
$f\theta\in L^1([0,\infty))\cap L^p([0,\infty))$ for some $p>\frac{1}{\alpha-1}$,
and $\beta+1>\alpha.$ Then, any continuous solution to (\ref{eqadn}) is asymptotically
($\mathcal{E},p$)-stable (resp. globally $\mathcal{E}$-stable in the large).
\end{theorem}
\noindent{\bf Proof}. Suppose that (H2) (resp. (H1)) is satisfied. By Theorem \ref{sim2}
we only need to see that the solution $Y$ of equation (\ref{n42}) is asymptotically 
($\mathcal{E},p$)-stable (resp. globally $\mathcal{E}$-stable in the large). Towards
this end, we invoke Lemma \ref{solueq} and \cite{FLM} (Lemma 2.4) to establish
\begin{equation*}
\begin{split}
Y(t)&=\xi_t+A\int_0^t (t-s)^{\beta-1}E_{\beta,\beta}(A(t-s)^\beta)\xi_sds\\
&\quad +\int_0^t (t-s)^{\alpha-1}E_{\beta,\alpha}(A(t-s)^\beta)f(s)d\theta_s\\
&=I_1(t)+I_2(t)+I_3(t),\quad t\ge 0.
\end{split}
\end{equation*}
Thus, considering Proposition \ref{ejemcomb} and \cite{FLM} (proof of Proposition 4.1) we only
need to show that, given $\varepsilon>0,$ $||I_3||_{\infty,[0,\infty)}<\varepsilon$ if
$||\xi||_{\mathcal{X}}+||f\theta||_{L^1([0,\infty))}+||f\theta||_{L^p([0,\infty))}$ $+
||\dot f\theta||_{L^1([0,\infty))}$ is small enough. For this purpose, we observe that 
(\ref{derimit}) and \cite{FLM} (Lemma 2.7) imply
\begin{eqnarray}\nonumber
I_3(t)&=&\int_0^t(t-s)^{\alpha-2}E_{\beta,\alpha-1}(A(t-s)^\beta)\theta_sf(s)ds\\ \nonumber
&&\quad -\int_0^t (t-s)^{\alpha-1}E_{\beta,\alpha}(A(t-s)^\beta)\theta_s\dot f(s)ds\\ \label{13.0}
&=& I_{3,1}(t)+I_{3,2}(t),\quad t\ge 0.
\end{eqnarray}
For $I_{3,1}$ we have, from (\ref{estiml}) and $q^{-1}=1-p^{-1},$
\begin{eqnarray}\nonumber
|I_{3,1}(t)| &\le & \int_0^{1\wedge t} s^{\alpha-2} |E_{\beta,\alpha-1}(As^\beta)| |\theta_{t-s}f(t-s)|ds
+\int_{1\wedge t}^t s^{\alpha-2} |E_{\beta,\alpha-1}(As^\beta)| |\theta_{t-s}f(t-s)|ds\\ \nonumber
&\le &  C_{\beta,\alpha-1}\left( \int_0^1s^{q(\alpha-2)}ds \right)^{1/q}
\left( \int_0^{1\wedge t} |\theta_{t-s} f(t-s)|^pds\right)^{1/p} +C_{\beta,\alpha-1}\int_0^t |\theta_{t-s}f(t-s)|ds\\
\label{13.4} &\le & C\left( ||\theta f||_{L^p([0,\infty)}+ ||\theta f||_{L^1([0,\infty)}  \right),\quad t\ge 0.
\end{eqnarray}

Finally, using (\ref{estiml}) again and the fact that $\beta+1>\alpha,$
\begin{equation*}
\begin{split}
I_{3,2}(t) &\le  \int_0^{1\wedge t} s^{\alpha-1} |E_{\beta,\alpha}(As^\beta)| |\theta_{t-s}\dot f(t-s)|ds+
\int_{1\wedge t}^t s^{\alpha-1} |E_{\beta,\alpha}(As^\beta)| |\theta_{t-s}\dot f(t-s)|ds\\
& \le C_{\beta,\alpha}\int_0^t |\theta_{t-s}\dot f(t-s)|ds+\frac{C_{\beta,\alpha}}{|A|}\int_0^t |\theta_{t-s}\dot f(t-s)|ds\\
&\le C\int_0^\infty |\theta_{t-s}\dot f(t-s)|ds,\quad t\ge 0.
\end{split}
\end{equation*}
Hence (\ref{13.0}) and (\ref{13.4}) yield that the proof is complete.

\hfill $\Box$

Observe that, in the previous proof, the inequality 
$$I_{3,2}(t)\le C\int_0^\infty |\theta_{t-s}\dot f(t-s)|ds,\quad t\ge 0,$$ 
is still true for $\beta+1\ge \alpha,$ wich is used in the proof of Theorem \ref{pro4.8} below.

\medskip

\subsection{Stochastic integral equations with additive noise}

In the remaining of this paper we suppose that all the introduced random variables are defined on a 
complete probability space $(\Omega,\mathcal{F},P).$

\begin{remark}\label{remarka}

Note that, in equation (\ref{eqadn}), we can consider a random variable $A:\Omega\to (-\infty,0),$
stochastic processes $\xi,$ $\theta$ and $f,$ and a random field $h$ such that for almost all $\omega,$
$A(\omega),$ $\xi_\cdot(\omega), \theta_\cdot(\omega), f(\omega,\cdot)$  and $h(\omega,\cdot)$
satisfy the hypotheses of Theorem \ref{sim3} (or Theorem \ref{sim2}), then we can analyze stability for
equation (\ref{eqadn}) $\omega$ by $\omega$ (i.e., with probability one). An example for the process $\theta$
is fractional Brownian motion $B^H$ with Hurst parameter $H\in (0,1)$. Fractional Brownian motion is a centered Gaussian
process with covariance 
$$R_{H}(s,t)=\mathbf{E}(B_s^HB_t^H)=\frac{1}{2}\big{(}s^{2H}+t^{2H}-|t-s|^{2H} 
\big{)},\quad s,t\ge 0.$$
It is well-known that $B^H$ has $\gamma$-H\"older continuous paths on compact sets,  for any exponent $\gamma<H$
due to Kolmogorov continuity theorem (see Decreusefond and \"Ust\"unel \cite{DU}).
\end{remark}

The last remark motivate the following:

\begin{definition} A continuous solution $X$ to equation (\ref{eqadn}) is said to be 
{\rm globally $\mathcal{E}$-stable in the mean} if $\mathbf{E} |X(t)|\to 0$ as $t\to\infty$
for any proces $\xi\in\mathcal{E}.$
\end{definition}

An immediate consequence of the proof of Theorem \ref{teolicon}, we can state the following extension of 
Theorem \ref{sim2}.

\begin{theorem}\label{asterisco} Let $h$ satisfy (H1), $A<0$, $\mathcal{E}$ a family of continuous processes
and $f,\theta$ as in Remark \ref{remarka} such that the solution to equation (\ref{n42}) is stable in
the mean. Then, any continuous solution to equation (\ref{eqadn}) is also $\mathcal{E}$-stable in the
mean.
\end{theorem}

\noindent{\bf Remark}. In \cite{FLM} (Theorem 4.3) we can find examples of families of processes for which the solution
of (\ref{n42}) is $\mathcal{E}$-stable in the mean.

Other definition motivated by Remark \ref{remarka} is the following:

\begin{definition} Let $\mathcal{E}\subset\mathcal{X}$ be a family of continuous functions. We say that
a continuous process $\xi$ belongs to $\mathcal{E}$ in the mean ($\xi\in\mathcal{E}_m$
for short) if $\mathbf{E}(|\xi|)\in\mathcal{E}$.
\end{definition}

Now we consider the stochastic integral equation
\begin{equation}\label{13.6}
X(t)=\xi_t+\frac{1}{\Gamma(\beta)}\int_0^t (t-s)^{\beta-1} [AX(s)+h(X(s))]ds+
\frac{1}{\Gamma(\beta+1)}\int_0^t (t-s)^{\beta}f(s)dB_s^\gamma, \quad t\ge 0.
\end{equation}
Here, in order to finish the paper, $A,h,\beta,\gamma$ and $f$ are as in equation (\ref{eqadn}) such that
$\beta+\gamma>1,$ and $\xi$ is a continuous stochastic process. We remark that we interprete equation
(\ref{13.6}) path by path (i.e. $\omega$ by $\omega$).

The following definition is also inpired by Remark \ref{remarka}.

\begin{definition}
Let $\mathcal{E}\subset\mathcal{X}$ be a family of continuous functions. We say that a continuous solution
to equation (\ref{13.6}) is ($\mathcal{E},p$)-stable in the mean if for a given $\varepsilon>0$ there is $\delta>0$
such that $||\mathbf{E}|X|||_{\infty,[0,\infty)}<\varepsilon$ for any $\xi\in\mathcal{E}_m$ such that

\begin{equation*}
\big{|}\big{|}\mathbf{E}|\xi|\big{|}\big{|}_{\mathcal{X}}+\big{|}\big{|}f(\cdot)\cdot^\gamma\big{|}\big{|}_{L^1([0,\infty))} 
+\big{|}\big{|}f(\cdot)\cdot^\gamma\big{|}\big{|}_{L^p([t_0,\infty))}
+\big{|}\big{|}\dot{f}(\cdot)\cdot^\gamma\big{|}\big{|}_{L^1([0,\infty)} <\delta.
\end{equation*}
\end{definition}   

\noindent{\bf Remark}. In this definition, if $\xi=\sum_{i=1}^n\xi^{(i)},$ with $\xi^{(i)}\in\mathcal{E}_m,$ then
we set $||\xi||_{\mathcal{X}}=\sum_{i=1}^n ||\xi^{(i)}||_{\mathcal{X}}.$

\begin{theorem}\label{pro4.8}
Let (H2) be true, $\xi$ as in Proposition \ref{ejemcomb}, $p>\frac{1}{\beta}$ and $f\in C^1((0,\infty))$ a positive
function with negative derivative such that $\left( r\mapsto r^\gamma |\dot{f}(r)|\right)\in L^1 ([0,\infty))$
and $\left( r\mapsto r^\gamma {f}(r)\right)\in L^1 ([0,\infty))\cap L^p([0,\infty)).$ Moreover, let $h$ be a 
non-decreasing and locally Lipschitz function, wich is concave on $\mathbb{R}_+$ and convex on $\mathbb{R}_-
\cup \{ 0\}.$ Then, the solution to equation (\ref{13.6}) is ($\tilde{\mathcal{E}},p$)-stable in the mean, where $\xi\in\tilde{\mathcal{E}}$ if and only if 
$\xi=\xi^{(1)}-\xi^{(2)}$ with $\xi^{1},\xi^{2}$ two non-negative, 
non-decreasing and continuous processes in
$\mathcal{E}_m.$
\end{theorem}
\noindent{\bf Proof}. Let $X$ be the continuous solution to equation (\ref{13.6}). Then Lemma \ref{solueq}
implies
\begin{equation*} 
\begin{split}
X(t)&=\xi_tE_{\beta,1}(At^\beta)+A\int_0^t(t-s)^{\beta-1}E_{\beta,\beta}(A(t-s)^\beta)(\xi_s-\xi_t)ds+
\int_0^t (t-s)^{\beta-1}E_{\beta,\beta}
(A(t-s)^\beta)h(X(s))ds\\ &\quad+\int_0^t (t-s)^{\beta}E_{\beta,\beta+1}(A(t-s)^\beta)f(s)dB_s^\gamma\\
&\le \xi_t^{(1)}E_{\beta,1}(At^\beta)+A\int_0^t (t-s)^{\beta-1}E_{\beta,\beta} (A(t-s)^\beta)(\xi_s^{(1)}-
\xi_t^{(1)})ds\\ &\quad+\int_0^t (t-s)^{\beta-1}E_{\beta,\beta}
(A(t-s)^\beta)h(X(s))ds+
\int_0^t (t-s)^{\beta-1}E_{\beta,\beta}(A(t-s)^\beta)|B_s^\gamma|f(s)ds\\
&\quad -\int_0^t (t-s)^{\beta}E_{\beta,\beta+1}(A(t-s)^\beta)\dot f(s)|B_s^\gamma|ds, \quad t\ge 0,
\end{split}
\end{equation*}
where the last inequality follows from the facts that $0\le \xi^{(1)}, \xi^{(2)}$ are two non-decreasing
processes, $f,(-\dot{f})\ge 0$ and from \cite{FLM} (Lemma 2.7). Therefore, we can state, by Lemma
\ref{pach}, that $X\le X^{(1)}$ where $X^{(1)}$ is the solution to
\begin{eqnarray}\nonumber
X^{(1)}(t)&=&\xi_t^{(1)}E_{\beta,1}(At^\beta)+A\int_0^t(t-s)^{\beta-1}E_{\beta,\beta}
(A(t-s)^\beta)(\xi_s^{(1)}-\xi_t^{(1)})ds
\\ \nonumber &&+\int_0^t (t-s)^{\beta-1}E_{\beta,\beta}(A(t-s)^\beta)h(X^{(1)}(s))ds+
\int_0^t (t-s)^{\beta-1}E_{\beta,\beta}(A(t-s)^\beta)|B_s^\gamma|f(s)ds\\  &&  \label{13.8.1}
 -\int_0^t (t-s)^{\beta}E_{\beta,\beta+1}(A(t-s)^\beta)\dot f(s)|B_s^\gamma|ds, \quad t\ge 0.
\end{eqnarray}
Observe that we also have $X^{(1)}(t)\ge 0$ due to $h(0)=0,$ Lemma \ref{pach} and
$$-X^{(1)}(t)\le \int_0^t (t-s)^{\beta-1}(A(t-s)^\beta)\hat{h}(-X^{(1)}(s))ds,\quad t\ge 0,$$
with $\hat{h}(x)=-h(-x),x\in\mathbb{R}.$ Proceeding similarly we have $-X(t)\le X^{(2)}(t),$
with $X^{(2)}(t)>0$ and
\begin{eqnarray}\nonumber
X^{(2)}(t)&=&\xi_t^{(2)}E_{\beta,1}(At^\beta)+A\int_0^t(t-s)^{\beta-1}E_{\beta,\beta}
(A(t-s)^\beta)(\xi_s^{(2)}-\xi_t^{(2)})ds
\\ \nonumber &&+\int_0^t (t-s)^{\beta-1}E_{\beta,\beta}(A(t-s)^\beta)\hat{h}(X^{(2)}(s))ds+
\int_0^t (t-s)^{\beta-1}E_{\beta,\beta}(A(t-s)^\beta)|B_s^\gamma|f(s)ds\\  &&  \label{13.8.2}
 -\int_0^t (t-s)^{\beta}E_{\beta,\beta+1}(A(t-s)^\beta)\dot f(s)|B_s^\gamma|ds, \quad t\ge 0.
\end{eqnarray}
In other words, we have
\begin{equation}\label{13.9.1}
\mathbf{E}\left(|X(t)|\right)\le \mathbf{E}\left(X^{(1)}(t)\right)+\mathbf{E}\left(X^{(2)}(t)\right),\quad t\ge 0.
\end{equation}

Finally, observe that (\ref{13.8.1}), (\ref{13.8.2}), the fact that $A$ is a negative number 
and Jensen inequality give, for $\theta_s=s^\gamma,$
\begin{equation*}
\begin{split}
\mathbf{E}\left(X^{(1)}(t)\right)
&\le\mathbf{E}(\xi_t^{(1)})E_{\beta,1}(At^\beta)+A\int_0^t(t-s)^{\beta-1}E_{\beta,\beta}
(A(t-s)^\beta)\mathbf{E}(\xi_s^{(1)}-\xi_t^{(1)})ds
\\  &\quad+\int_0^t (t-s)^{\beta-1}E_{\beta,\beta}(A(t-s)^\beta)h(\mathbf{E}[X^{(1)}(s)])ds
+\int_0^t (t-s)^{\beta}E_{\beta,\beta+1}(A(t-s)^\beta) f(s)d\theta_s, \quad t\ge 0,
\end{split}
\end{equation*}
and 
\begin{equation*}
\begin{split}
\mathbf{E}
\left(X^{(2)}(t)\right)&\le\mathbf{E}(\xi_t^{(2)})E_{\beta,1}(At^\beta)+A\int_0^t(t-s)^{\beta-1}E_{\beta,\beta}
(A(t-s)^\beta)\mathbf{E}(\xi_s^{(2)}-\xi_t^{(2)})ds
\\  &\quad+\int_0^t (t-s)^{\beta-1}E_{\beta,\beta}(A(t-s)^\beta)
\hat{h}(\mathbf{E}[X^{(2)}(s)])ds
+\int_0^t (t-s)^{\beta}E_{\beta,\beta+1}(A(t-s)^\beta) f(s)d\theta_s, \quad t\ge 0,
\end{split}
\end{equation*}

Hence by (\ref{13.9.1}), Lemma \ref{pach}, Hypothesis (H2), and the proofs of Proposition \ref{ejemcomb} and Theorem
\ref{sim3} we get that the result holds. Indeed, for $i=1,2,$
$$\mathbf{E}(X^{(i)}(t))\le u^{(i)}(t),\quad t\ge 0,$$
where $u^{(i)}$ is the unique solution to the equation
$$u^{(i)}(t)=\mathbf{E}(\xi^{(i)}_t)+\frac{1}{\Gamma(\beta)}\int_0^t (t-s)^{\beta-1}[A+C]u^{(i)}(s)ds+
\frac{1}{\Gamma(\beta+1)}\int_0^t (t-s)^\beta f(s)d\theta_s,\quad t\ge 0.$$
\hfill $\Box$

\begin{example}
{\rm A function $h$ that satisfies the conditions of Theorem \ref{pro4.8} is}
$$h(x) = \begin{cases} 1-e^{-Cx}, & \mbox{if } x\ge 0 \\ e^{Cx}-1, & \mbox{if } x<0, \end{cases}$$
{\rm where $C>0.$ Indeed, we have that}
$$h'(x) = \begin{cases} Ce^{-Cx}, & \mbox{if } x\ge 0 \\ Ce^{Cx}, & \mbox{if } x<0. \end{cases}$$
{\rm Thus, given $\varepsilon>0$ there is $\delta>0$ such that 

$$|h(x)|\le (C+\varepsilon)|x| \quad \mbox{for}\quad |x|\le\delta.$$}

\end{example}

\bigskip

\begin{example}{\rm Here we give a function that satisfies Assumption 2 on Definition \ref{coninic}.
Let $\xi_t=g(t)\sin{\frac{1}{t}},$ $t\ge 0.$
The function $g$ is bounded and satisfies $g(t)=\psi(t)c_0t^{3-\upsilon}+\varphi(t)\frac{c_1}{1+t},$
where $\psi,\varphi\in C^\infty(\mathbb{R}_+)$ are such that 
\[ \psi(t) = \left\{ \begin{array}{ll}
         1 & \mbox{if $t \in [0,1]$};\\
        0 & \mbox{if $t \ge  2$},\end{array} \right. \mbox{and} \quad \varphi(t) = \left\{ \begin{array}{ll}
         0 & \mbox{if $t \in [0,1]$};\\
        1 & \mbox{if $t \ge  2$}.\end{array} \right. \]         
Thus}
\end{example}
$$\xi'_t =g'(t)\sin{\frac{1}{t}}-g(t)t^{-2}\cos{\frac{1}{t}},\quad t\ge 0.$$
Now it is easy to verify our claim is true using straightforward calculations.

\medskip

\noindent{\bf Aknowledgements}: The authors thank Cinvestav-IPN and Universitat de Barcelona
for their hospitality and economical support.


\begin{thebibliography}{99}


\bibitem{AF} J.A.D.Appleby and A.Freeman:  {\it Exponential asymptotic stability of linear It\^o-Volterra equations with
damped stochastic perturbations}, Electron. J. Probab. {\bf 8}, (2003) 1-22. 

\bibitem{Ba} R.L.Bagley and R.A.Calico: {\it Fractional order state equations for the control
of viscoelastically damped structures,} J. Guid. Control Dyn.,  {\bf 14}, (1991) 304-311.

\bibitem{B} T.Q.Bao: {\it On the existence, uniqueness and stability of solutions of stochastic Volterra-Ito 
equation}, Vietnam Journal of Mathematics {\bf  32}, (2004) 389-397.


\bibitem{Da} S.Das: {\it Functional Fractional Calculus for Systems Identication and Controls,}
Springer Berlin Heidelberg New York, 1sth edition (2008).

\bibitem{DEOK} A.Demir, S.Erman, B.\"Ozg\"ur and E.Korkmaz: {\it Analysis of fractional partial differential
equations by Taylor series expansion}, Boundary Value Problems (2013).


\bibitem{DLC} W.Deng, C.Li, L.Changpin and J.L\"u: {\it Stability analysis of linear fractional differential system with multiple time delays}, Nonlinear Dynamics,  {\bf 48}, (2007) 409-416.


\bibitem{DU} L.Decreusefond and  A.S.\"Ust\"unel:  {\it Stochastic analysis of the fractional 
Brownian motion},  Potential Anal.  {\bf 10}, (1999) 177-214.

\bibitem{DNo}
R.M.Dudley and R.Norvai\u{s}a, \textit{An introduction to $p$-variation and 
Young Integrals}, Tech. Rep. 1, Maphysto, Centre for Mathematical Physics and Stochastics, 
University of Aarhus. Concentrated advanced course (1998).

\bibitem{FH} P.K.Friz and M.Hairer: {\it A Course on Rough Paths: With an Introduction to Regularity Structures},
Springer International Publishing Switzerland (2014).

\bibitem{FLM} A.Fiel, J.A.Le\'on and D.M\'arquez-Carreras: {\it 
Stability for some linear stochastic fractional systems}, Communications on Stochastic Analysis, Serials Publications,
{\bf 8},  (2014) 205-225. 

\bibitem{G} M.Gubinelli: {\it Controlling rough paths,} J. Funct. Anal. 
{\bf 216}, (2004) 86-140.

\bibitem{GG} I.Grigorenko and  E.Grigorenko: {\it Chaotic Dynamics of the Fractional
Lorenz System,} Physical Review Letters,  {\bf 91} (2003).

\bibitem{HLQK} T.T.Hartley, C.F.Lorenzo and H.K.Qammer: {\it Chaos in a fractional order Chua`s system}, 
IEEE Transactions on circuits and systems-I: Fundamental theory and applications,  {\bf 42}, (1995) 485-490.

\bibitem{Ho} O.Heaviside:{\it Electromagnetic Theory}, New York: Chelsea (1971).

\bibitem{H} R.Hilfer: {\it Applications of Fractional Calculus in Physics,} World Scientic, River
Edge, New Jersey (2000).

\bibitem{INK} M.Ichise, Y.Nagayanagi, T.Kojima: {\it An analog simulation of non-integer
order transfer functions for analysis of electrode processes,} J. Electroanal. Chem.,
 {\bf 33}, (1971) 253-265.

\bibitem{JJM} D.Junsheng, A.Jianye and X.Mingyu: {\it Solution of system of 
fractional differential equations by adomian decomposition method,} {\rm Appl. 
Math. J. Chinese Univ. Ser. B} {\bf 22 (1)}, (2007) {7-12}.

\bibitem{KST} A.A.Kilbas, H.M.Srivastava and J.J.Trujillo: 
{\it Theory and Applications of Fractional Differential Equations: Theory and Applications},
Elsevier B. V. (2006).

\bibitem{KBD} D.Kusnezov,  A.Bulgac and G.D.Dang: {\it Quantum Levy processes and fractional kinetics}, Phys. Rev. Lett., 
{\bf 82}, (1999) 1136-1139.

\bibitem{LV} V.Lakshmikantham and A.S.Vatsala: {\it Basic theory of fractional differential equations},
Nonlinear Analysis {\bf 69}, (2008) 2677-2682.

\bibitem{LT} J.A.Le\'on and S.Tindel: {\it Malliavin calculus for fractional delay equations},
{\rm J. Theoret. Probab.} {\bf 25 (3)}, (2012) {854-889}.

\bibitem{LLW} W.Li, M.Liu and K.Wang:  A generalization of It\^o's formula and the stability of stochastic
Volterra integral equations, {\em Journal of Applied Mathematics} Article ID 292740  (2012). 


\bibitem{LCP} Y.Li , Y.Q.Chen,  and I.Podlubny: {\em Stability of fractional-order
systems: Lyapunov direct method and generalized Mittag-Leffler stability,} 
Computer and Mathematics with Applications {\bf 59}, (2010) 1810-1821.

\bibitem{Lin} S.J.Lin: {\it Stochastic analysis of fractional Brownian motions}, Stochastics Stochastics Rep. {\bf 55}, 
(1995) 121-140.


\bibitem{Ly} T.Lyons: {\it T. Lyons: Differential equations driven by rough signals (I): An extension
of an inequality of L. C. Young}, Mathematical Research Letters {\bf 1}, (1994)
451-464.

\bibitem{MA} J.A.T.Machado and A.Azenha: {\it Fractional-order hybrid control of robot manipulators,} 
En: Systems, Man, and Cybernetics, 1998. 1998 IEEE International Conference on,  {\bf 1},  (1998) 788-793.

\bibitem{M}  D.Matignon: {\it Stability results for fractional differential equations with applications to
control processing}, {\em In Proc IMACS, IEEE-SMC}, (1996) 963-968.  

\bibitem{MM-L-FA} R.Mart\'{\i}nez-Mart\'{i}nez, J.A.Le\'on and G.Fern\'andez-Anaya: 
{\it Asymptotc stability of fractional order nonlinear systems via Lyapunov like conditions},
{\rm Preprint} (2013).


\bibitem{MS} K.S.Miller and S.G.Samko: {\it Completely monotonic functions,} Integr. Transf. and 
Spec. Funct.,   {\bf 12}, (2001) 389-402.

\bibitem{MO} S.Momania and Z.Odibat: {\it A novel method for nonlinear fractional 
partial differential equations: Combination of DTM and generalized Taylor's formula}, 
Journal of Computational and Applied Mathematics {\bf 220}, 
(2008) 85-95.

\bibitem{MO2} S.Momania and Z.Odibat: {\it Numerical methods for nonlinear partial 
differential equations of fractional order}, Applied Mathematical Modelling {\bf 32} (2008) 28-39.

\bibitem{DN} D.Nguyen: {\em Asymptotic behavior of linear fractional stochastic differential equations with 
time-varying delays,}  Commun Nonlinear Sci. Numer. Simul. {\bf 19}, (2014) 1-7.

\bibitem{N1} D.Nualart: {\it Stochastic integration with respect to fractional
Brownian motion and applications}, {\em In Stochastic Models, Contemporary Mathematics} {\bf 336},
(2003) 3-39.




\bibitem{NR}  D.Nualart and A.R\u{a}\c{s}canu: {\it Differential equations driven by 
fractional Brownian motion}, {\em Collect Math.} {\bf 53},  (2002) 55-81.

\bibitem{Pa} B.G.Pachpatte: {\it Inequalities for Differential and Integral Equations}, 
Academic Press Limited,  {\bf 197}, Marathwada University, Aurangabad, India (1998). 

\bibitem{Pod1} I.Podlubny: {\it Fractional Differential Equations}. {\rm 9th Edition, 
Academic Press}, San Diego (1999).

\bibitem{QT} L.Quer-Sardanyons and S.Tindel: {\it Pathwise definition 
of second-order SDES},  Stochastic Processes and their Applications
{\bf 122}, (2012) 466-497.

\bibitem{RSES}   A.G.Radwan, A.M.Soliman,  A.S.Elwakil and A.Sedeek: {\it On the stability of linear systems 
with fractional-order elements}, Chaos, Solitons and Fractals {\bf 40}, (2009) 2317-2328. 


\bibitem{RV} F.Russo and P.Vallois: \textit{Forward, backward and symmetric 
stochastic integration}, Probab. Theory Relat. Fields {\bf 97}, (1993) 403-421.

\bibitem{S} R.Schneider: {\it Completely monotone generalized Mittag-Leffler 
functions}, Expo. Math {\bf 14}, (1996) 3-16.

\bibitem{WWL} X.-J.Wen, Z.-M.Wu and J.-G-Lu: {\it Stability analysis of a class 
of nonlinear fractional-order systems,} {\rm IEEE Transactions on circuits and 
systems - II: Express Briefs.} {\bf 55 (11)}, (2008) 1178-1182.

\bibitem{YL} J.Yan and C.Li: {\it On chaos synchronization of fractional differential equations}, Chaos Solutions \& Fractals (2005).

\bibitem{YZ} Z.Yan and H.Zhang: {\it Asymptotic stability of fractional impulsive neutral stochastic
 integro-differential equations with state-dependent delay,} Electronic Journal of Differential Equations, Vol. 2013, 
 (2013) 1-29.

\bibitem{YLS} Y.Yu,  H.Li, and Y.Su: {\it The Synchronization of Three Chaotic Fractional-order Lorenz 
Systems with Bidirectional Coupling}, Journal of Physics: Conference Series (2007).

\bibitem{Y} L.C.Young: {\it An inequality of the H\"older type, connected with Stieltjes 
integration,} Acta Math. {\bf 67}, (1936) 251-282. 

\bibitem{ZC} F.Zhang and C.Li: {\it Stability Analysis of Fractional Differential Systems with Order Lying in $(1,2),$}
Advances in Difference Equations (2011).

\bibitem{ZL} C.Zhang, W.Li and K.Wang:{\it Stability and boundedness of stochastic Volterra integrodifferential 
equations with infinite delay},  Journal of Applied Mathematics, Article ID 320832  (2013).

\bibitem{ZZ} B.Zhang and J.Zhang:{\it Conditional stability of stochastic Volterra equations with anticipating kernel}, 
Journal of Mathematical Research \& Exposition  {\bf 22(2)} (2002).

\bibitem{Z} M.Z\"ahle: {\it Integration with respect to fractal functions
and stochastic calculus I,} Probab. Theory Relat. Fields {\bf 111}, (1998) 333-374.

\bibitem{ZYC} C.Zeng, Q.Yang and Y.Q.Chen: {\it Lyapunov techniques for stochastic differential equations
driven by fractional Brownian motion,}  Abstract and Applied Analysis, Article ID 292653  (2014).

\end{thebibliography}
\end{document}